\documentclass[12pt,a4paper]{amsart}

\usepackage{amssymb}
\makeatletter
\renewcommand{\@begintheorem}[2]{
\rm \trivlist \item [\hskip \labelsep {\bf #2\ \ #1.}]
                                }
\makeatother

\makeatletter

\DeclareFontFamily{U}{cyr}{}
\DeclareFontShape{U}{cyr}{m}{n}{
  <5> wncyr5 <6> wncyr6 <7> wncyr7 <8> wncyr8 <9> wncyr9 <10->
wncyr10}{}
\DeclareMathAlphabet{\mathcyr}{U}{cyr}{m}{n}

\input cyracc.def
\input epsf

\newcommand{\ts}{\vspace{\baselineskip}\noindent{\bf Proof.}$\;\;$}
\newcommand{\ZZ}{{\bf Z}}
\newcommand{\QQ}{{\bf Q}}

\newcommand{\CC}{{\bf C}}

\newcommand{\PP}{{\bf P}}

\newcommand{\BB}{{\mathbb B}}

\newcommand{\rmd}{\mbox{d}}

\newcommand{\bes}{\begin{equation*}}
\newcommand{\ees}{\end{equation*}}

\def\DynkinEEE#1#2#3#4#5#6#7
   {\begin{picture}(124, 100)%
   \put(5, 80){\circle{4}}%
   \put(7, 80){\line(1, 0){28}}%
   \put(37, 80){\circle{4}}%
   \put(39, 80){\line(1, 0){28}}%
   \put(69, 80){\circle{4}}%
   \put(71, 80){\line(1, 0){28}}%
   \put(101, 80){\circle{4}}%
   \put(103, 80){\line(1, 0){28}}%
   \put(133, 80){\circle{4}}%
   \put(69, 78){\line(0, -11){28}}%
   \put(69, 48){\circle{4}}%
   \put(69, 46){\line(0, -11){28}}%
   \put(69, 16){\circle{4}}%
   \put(5, 86){\makebox(0, 0)[b]{$#1$}}%
   \put(37, 86){\makebox(0, 0)[b]{$#3$}}%
   \put(69, 86){\makebox(0, 0)[b]{$#4$}}%
   \put(101, 86){\makebox(0, 0)[b]{$#5$}}%
   \put(133, 86){\makebox(0, 0)[b]{$#6$}}%
   \put(75, 48){\makebox(0, 0)[lb]{$#2$}}%
   \put(75, 16){\makebox(0, 0)[lb]{$#7$}}%
\end{picture}}

\title{The Picard-Fuchs equation of a family of Calabi-Yau threefolds without maximal unipotent monodromy}
\author{Alice Garbagnati}
\author{Bert van Geemen}
\address{Dipartimento di Matematica, Universit\`a di Milano,
Via Saldini 50, I-20133 Milano, Italia}
\email{alice.garbagnati@unimi.it}
\email{lambertus.vangeemen@unimi.it}

\begin{document}

\begin{abstract}
Recently J.C.\ Rohde constructed families of Calabi-Yau
threefolds parametrised by Shimura varieties. The points corresponding to threefolds with CM are dense in the Shimura variety and, moreover, the families do not have boundary points with maximal unipotent monodromy. Both aspects are of interest for Mirror Symmetry. In this paper we discuss one of Rohde's examples in detail and  we explicitly give the Picard-Fuchs equation for this one dimensional family.
\end{abstract}

\maketitle

In this note we work out an example of J.C.\ Rohde of a one dimensional family of Calabi-Yau threefolds $X_\lambda$ (with $h^{2,1}(X_\lambda)=1$), 
parametrised by a Shimura variety, such that the associated Picard-Fuchs equation has no maximal unipotent monodromy. Actually, as already pointed out by Rohde, the Picard-Fuchs equation of degree four reduces to two differential equations of degree two. In section \ref{hge} we give these two
equations, which are hypergeometric differential equations, explicitly.

This is of some interest for Mirror Symmetry, which suggests that
a family like the $X_\lambda$ should be the Mirror of another family of Calabi-Yau threefolds $Y_\mu$. A particular solution of the Picard-Fuchs equation, chosen using the maximal unipotent monodromy, should have a Taylor expansion whose $d$-th coefficient is related to the number of rational curves of 
degree $d$ on a general $Y_\mu$ (more precisely, it is related to a certain Gromow-Witten invariant of $Y_\mu$).
In the absence of maximal unipotent monodromy, the recipe for identifying the solution of the Picard-Fuchs equation should be modified. It might however be the case that the family $X_\lambda$ is not the Mirror of any family of Calabi-Yau's, which would also be of some interest.

We follow the approach indicated by Rohde in \cite{Rohde-maximal}. The main point is a good understanding of a family of K3 surfaces which is used in the construction of the $X_\lambda$'s. This leads to an explicit description of 
$H^3(X_\lambda,\QQ)$ in terms of $H^1(C_\lambda,\QQ)$ for a genus two curve
$C_\lambda$. The only thing to do then is to recall a classical result on the Picard-Fuchs equations of the $C_\lambda$. 

In section \ref{cm} we briefly comment on the 
Calabi-Yau threefolds $X_\lambda$ with complex multiplication (CM), in particular we observe that these are related to elliptic curves with CM. It is conjectured that Calabi-Yau threefolds with CM are related to Rational Conformal Field Theories (RCFT) and these should be easier to understand then the more general
Conformal Field Theories (cf.\ \cite{GV}).

\section{The one parameter family of Calabi-Yau threefolds.}

\subsection{Rohde's construction}\label{rohde}
The construction of Rohde  (\cite{Rohde-maximal}, section 4) 
starts with a K3 surface $S$ with an automorphism $\alpha$ which acts
by $\overline{\xi}$ on $H^{2,0}(S)\cong\CC$, where $\xi=e^{2\pi i/3}$ is a primitive cube of unity.  
To get from such a K3 surface to a Calabi-Yau threefold, one assumes that
the fixed point locus of $\alpha$ consists of $k$ disjoint smooth rational curves and $k+3$ isolated fixed points. 
These K3 surfaces are parametrised by the $r$-ball $\BB_r\cong SU(1,r)/S(U(1)\times U(r))$, where 
$r=6-k$, see \cite{AS}, \cite{Rohde-maximal} and also
\cite{DK} and the references given there for more results of this kind.

Let $E\cong \CC/\ZZ+\ZZ\xi$ be the elliptic curve with $j$-invariant zero, so $E$ is isomorphic to the Fermat cubic curve and has Weierstrass equation $y^2=x^3+1$. We define an automorphism of $E$ by
$$
\alpha_E\,:\,E\longrightarrow E,\qquad (x,y)\,\longmapsto\, (\xi x,y),
\qquad(E\,:\,y^2=x^3+1)
$$ 
it acts as multiplication by $\xi$ on $H^{1,0}(E)=\CC\rmd x/y$.
Then the product $S\times E$ has the automorphism $(\alpha,\alpha_E)$
which is trivial on $H^{2,0}(S)\otimes H^{1,0}(E)$. Rohde shows that  the quotient 
$(S\times E)/(\alpha,\alpha_E)$ is birationally isomorphic to a Calabi-Yau threefold ${X}$ with 
$$
h^{2,1}({X})\,=\,r\,=\,6-k,\qquad  
h^{1,1}({X})\,=\,18+11k.
$$
These Calabi-Yau varieties are then parametrised by the ball $\BB_r$ and their moduli space is a quotient of the ball by an arithmetic subgroup of $SU(1,r)$.

\subsection{The case $r=1$.}\label{r=1}
We consider the  family of K3 surfaces $S_f$,
which have an elliptic fibration $\pi:S_f\rightarrow \PP^1_t$ with Weierstrass model
$$
Y^2\,=\,X^3\,+\,f(t)^2,\qquad f:=gh^2,\quad \mbox{deg}(g)=\mbox{deg}(h)=2,
$$
and we assume that $g(t)h(t)$ has 4 distinct zeroes in 
$\CC\subset \PP^1_t$.
Each surface $S_f$ has an automorphism $\alpha_f$ of order three defined
by
$$
\alpha_f\,:\,S_f\,\longrightarrow\,S_f,\qquad
(X,Y,t)\,\longmapsto\, (\overline{\xi}X,Y,t)
$$
which acts as $\overline{\xi}$ on 
$H^{2,0}(S_f)\,=\,\CC \rmd t\wedge \rmd X/Y$.
This automorphism fixes the section at infinity $s_\infty$, a smooth rational curve, and maps each fiber of $\pi$ into itself. 
The elliptic fibration has two other sections, $s_\pm(t):=(0,\pm f(t),t)$, of order three in the Mordell-Weil group, which are fixed under $\alpha_f$. There are four singular fibers over the zeroes of $f$. 
As the $j$-invariant of the smooth fibers is $j=0$ and the discriminant of the Weierstrass model has zeroes of order four resp.\ eight in the
zeroes of $g$, $h$ resp., Table 4.1 in \cite{Silverman} determines their type.
The two fibers over $g=0$ are of type $IV$. 
Such a fiber consists of three smooth rational curves meeting in one point, moreover each of these curves meets one of the sections.
As the fixed point locus of $\alpha_f$ is smooth
and the sections are pointwise fixed, 
the components of such a fiber are mapped into themselves by $\alpha_f$ 
and $\alpha_f$ induces a non-trivial automorphism on each of these. 
Thus their point of intersection is a fixed point of $\alpha_f$ in $S$.
The two fibers over $h=0$ are of type $IV^*$, so they consist of 7 smooth rational curves with intersection graph (which is the affine Dynkin diagram of type $\tilde{E}_6$):
$$
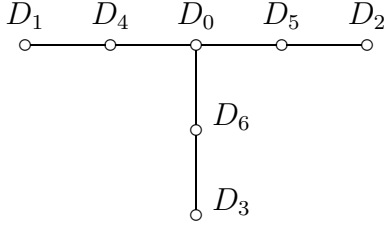
\DynkinEEE{D_1}{D_6}{D_4}{D_0}{D_5}{D_2}{D_3}
$$
Each of the three sections meets one of the curves $D_1,D_2,D_3$, so these curves, and hence all seven curves in the $IV^*$ fiber, are mapped into themselves under $\alpha_f$. 
Moreover, $\alpha_f$ is non-trivial on $D_1,D_2,D_3$, and thus has two fixed points on $D_1,D_2,D_3$ where these meet a section and another component of the fiber. 
The three points $D_4\cap D_0,D_5\cap D_0$ and $D_6\cap D_0$ on $D_0$ are fixed by $\alpha_f$, hence $D_0$ is pointwise fixed by $\alpha_f$.
Then $\alpha_f$ is non-trivial on $D_4,D_5,D_6$ and the three intersection points 
$D_i\cap D_{i+3}$, $i=1,2,3$ are isolated fixed points.
Thus $\alpha_f$ has one fixed curve and three fixed points in each $IV^*$ fiber. Therefore the fixed point locus of $\alpha_f$ consists of
$3+2\cdot 1=5$ smooth rational curves and $2\cdot 1+2\cdot 3=8$ points, so $k=5$ and the number of moduli is $r=1$.

The N\'eron Severi group of $S_f$ contains the classes of a fiber, $s_\infty$, two of the three rational curves in each type $IV$ fiber and six of the seven rational curves in each type $IV^*$ fiber.
These classes span a lattice of rank $1+1+2\cdot 2+2\cdot 6=18$ and thus
the transcendental lattice $NS(S_f)^\perp\;(\subset H^2(S_f,\ZZ))$ has rank at most $4$. Since all the $S_f$ have an automorphism of order three and one modulus, 
the N\'eron Severi group of the general $S_f$ has rank $18$ and the eigenvalues of $\alpha_f$ on the orthogonal complement of this rank $18$ lattice must be $\xi,\overline{\xi}$, each with multiplicity 2.
The decomposition of $H^2(S_f,\QQ)$ into eigenspaces for $\alpha_f$ is thus as follows:
$$
H^2(S_f,\QQ)\,=\,H^2(S_f,\QQ)_1\,\oplus T_f,
\qquad T_f\,:=\,(H^2(S_f,\QQ)_1)^\perp,
\quad \dim H^2(S_f,\QQ)_1\,=\,18,
$$
and the complexification of  $T_f$ decomposes into $\alpha_f$-eigenspaces:
$$
T_f\otimes_\QQ\CC\,=\,T_{a,\overline{\xi}}\,\oplus\,\,T_{a,\xi}\,=
T^{2,0}_{a,\overline{\xi}}\oplus T^{1,1}_{a,\overline{\xi}}\oplus 
T^{1,1}_{a,\xi}\oplus T^{0,2}_{a,\xi},
$$ 
each of these four spaces is one dimensional.

\subsection{Remark}\label{g=3,4}
A general Weierstrass model as above with $f=t(t-1)(t-a_1)(t-a_2)$ defines for general $a_1,a_2$ an elliptic fibration on a K3 surface with 4 fibers of type $IV$ and one of type $IV^*$ (over $t=\infty$).
Thus the fixed point locus of $\alpha_f$ consists of $4$ smooth rational curves and $7$ points, hence $k=4$ and the number of moduli is $r=2$.

A Weierstrass model as above with $f$ a general polynomial of degree six in $t$ defines an elliptic fibration on a K3 surface with 6 fibers of type $IV$. 
Thus the fixed point locus of $\alpha_f$ consists of $3$ smooth rational curves and $6$ points, hence $k=3$ and the number of moduli is $r=3$. 

In case the two zeroes of $g$ coincide, so $f$ has three double zeroes,
the elliptic fibration has three fibers of type $IV^*$ and $r=0$. This surface is birational to the quotient surface 
$(E\times E)/(\alpha_E,\alpha_E^{-1})$, cf.\ the next section and 
\cite{SI}, Lemma 5.1.

\section {The Picard-Fuchs equation}

\subsection{The middle cohomology}
The construction of Rohde (cf.\ section \ref{rohde}) applied to a K3 surface $S_f$ from section \ref{r=1}
defines a Calabi-Yau threefold ${X}_f$.
To determine the Picard-Fuchs equation of this family, we first observe that the cohomology group $H^3({X}_f,\QQ)$ is obtained by an elementary construction from the cohomology group $H^1(C_f,\QQ)$ of a curve $C_f$.

This curve $C_f$ is the smooth, projective, genus two curve defined by the equation
$$
C_f\,:\quad v^3=f(t),\qquad f=gh^2,
$$
with $g,h$ as in section \ref{r=1}. 
This equation exhibits $C_f$ as a cyclic degree three cover of $\PP^1_t$. The covering automorphism is 
$$
\beta_f\,:\,C_f\longrightarrow\,C_f, \qquad
(t,v)\,\longmapsto\, (t,\xi v).
$$
Let $H^{p,q}(C_f)_\xi, H^{p,q}(C_f)_{\overline{\xi}}$ be the subspaces of $H^{p,q}(C_f)$ on which $\beta_f$ acts as $\xi,\overline{\xi}$ respectively. Then we have the following decomposition of $H^1(C_f,\CC)$ into in four one-dimensional eigenspaces:
$$
H^1(C_f,\CC)\,=\,H^{1,0}(C_f)_\xi\oplus H^{1,0}(C_f)_{\overline{\xi}}
\oplus  H^{0,1}(C_f)_\xi\oplus H^{0,1}(C_f)_{\overline{\xi}}.
$$
Note that $\overline{H^{p,q}(C_f)_\xi}=H^{q,p}(C_f)_{\overline{\xi}}$.
To be explicit, one has
$$
H^{1,0}(C_f)_{\xi}\,=\,\CC g\rmd t/v^2,\qquad
H^{1,0}(C_f)_{\overline{\xi}}\,=\,\CC\rmd t/v,.
$$

\subsection{Proposition}\label{cohX}
Let ${X}_f$ be the Calabi-Yau threefold which is birational to $(S_f\times E)/(\alpha_f,\alpha_E)$.
Then there is an isomorphism
$$
\phi\,:\,
H^1(C_f,\QQ) \,\stackrel{\cong}{\longrightarrow}\, H^3({X}_f,\QQ)
$$ 
such that 
$$
\phi(H^{1,0}({C}_f)_{\overline{\xi}})\,=\,H^{3,0}(X_f),\qquad
\phi(H^{0,1}({C}_f)_{\overline{\xi}})\,=\,H^{2,1}(X_f).
$$

\ts
All smooth fibers of the elliptic fibration $\pi:S_f\rightarrow \PP^1_t$ are isomorphic to $E$. 
Pulling back this fibration along the degree three cover 
$C_f\rightarrow \PP^1_t$ 
one obtains an elliptic fibration over $C_f$ with Weierstrass model $Y^2=X^3+v^6$. This fibration over $C_f$ is birational to the product $C_f\times E$ (via $(X,Y,t)\mapsto ((t,v),(x,y)):=((t,v),(v^{-2}X,v^{-3}Y))$). 

Moreover, the rational map
$$
C_f\times E\,\longrightarrow\,S_f,\qquad
((t,v),(x,y))\longmapsto (X,Y,t)=(v^2x,v^3y,t)
$$
identifies $S_f$ with the minimal model of the quotient surface
$(C_f\times E)/(\beta_f,\alpha_E)$.
The automorphism $\alpha_f$ of $S_f$ is induced by the automorphism
$(1,\alpha^{-1}_E)$ of $C_f\times E$.

As $H^{2,0}(S_f)\subset T_f^{2,0}$, one obtains an isomorphism
$$
T_f\,\cong\,\left(H^1(C_f,\QQ)\otimes H^1(E,\QQ)\right)^{(\beta_f,\alpha_E)}\qquad
(\subset\;H^2(S_f,\QQ)).
$$

Rohde showed that (cf.\ the proof of Proposition 4.5 in \cite{Rohde-maximal}):
$$
H^3({X}_f,\QQ)\,\cong\,(H^2(S_f,\QQ)\otimes H^1(E,\QQ))^{(\alpha_f,\alpha_E)}\,\cong\, 
(T_f\otimes H^1(E,\QQ))^{(\alpha_f,\alpha_E)},
$$
where one uses that $\alpha_E$ has no eigenvalue $1$ in $H^1(E,\QQ)$.
Combining this with the description of $T_f$, 
this leads to
$$
H^3({X}_f,\QQ)\,\cong\,
\left( H^1(C_f,\QQ)\otimes H^1(E,\QQ)\otimes H^1(E,\QQ)
\right)^G,\qquad 
G:=\langle (1,\alpha_E,\alpha_E^{-1}),(\beta_f,\alpha_E,1)\rangle.
$$
In particular, the variation of the Hodge structures on $\{H^3(X_f,\QQ)\}_f$ comes from the one on $\{H^1(C_f,\QQ)\}_f$.

Now we determine the $G$-invariants by first considering the last two tensor factors. Let
$$
T_E\,:=\,
\left(H^1(E,\QQ)\otimes H^1(E,\QQ)\right)^{(\alpha_E,\alpha_E^{-1})},
$$
this Hodge substructure has dimension two and the automorphism $(\alpha_E,1)$ induces 
an automorphism $\beta_E$ of $T_E$ of order three. 
We have the eigenspace decomposition for $\beta_E$:
$$
T_E\otimes_\QQ\CC\,=\,T^{2,0}_{E,\xi} \,\oplus\, T^{0,2}_{E,\overline{\xi}},\qquad
\dim T^{2,0}_{E,\xi}\,=\,\dim T^{0,2}_{E,\overline{\xi}}\,=\,1.
$$
Thus $H^3({X}_f,\QQ)$ is the Hodge substructure of $(\beta_f,\beta_E)$-invariants 
in $H^1(C_f,\QQ)\otimes T_E$.
In particular, 
$H^{3,0}({X}_f)=H^{1,0}(C_f)_{\overline{\xi}}\otimes T^{2,0}_{E,\xi}\cong H^{1,0}(C_f)_{\overline{\xi}}$ etc.
To get the map $\phi$, one observes that the $\QQ$-vector space 
$H^1(C_f,\QQ)\otimes T_E$ has the decomposition into 
$(\beta_f,\beta_E)$-stable subspaces:
$$
H^1(C_f,\QQ)\otimes T_E\,=\,
\left(H^1(C_f,\QQ)\otimes T_E\right)^{(\beta_f,\beta_E)}\,
\oplus V \,\cong\,H^{3}({X}_f,\QQ)\,\oplus\,V
$$
where $V\otimes\CC$ is the direct sum of the $\xi,\overline{\xi}$ eigenspaces of $(\beta_f,\beta_E)$. Fix a non-zero element $t\in T_E$. Then one defines $\phi$ to be the composition of 
$H^1(C_f,\QQ)\rightarrow H^1(C_f,\QQ)\otimes T_E$, $x\mapsto x\otimes t$
with the projection onto the summand $H^{3}({X}_f,\QQ)$.
\qed

\subsection{Remarks} Generalisations of tensoring with $H^1(E,\QQ)$ and taking invariants, as we used repeatedly in the proof of the proposition above, are considered in \cite{vG} and \cite{DK}, section 13.

The proposition can be easily generalised to the other cases considered in Remark \ref{g=3,4}.

Note that the proposition shows that $X_f$ is birationally isomorphic 
to $(C_f\times E\times E)/G$ and that the $G$-invariant holomorphic
three form $(\rmd t/v)\wedge (dx_1/y_1)\wedge(dx_2/y_2)$ descends to the holomorphic three form on $X_f$, where $(x_i,y_i)$ are the coordinates on the $i$-th copy of $E$. As the elliptic curve $E$ is fixed, the variation of the genus two curves $C_f$ determines the variation of the CY threefolds $X_f$.

\subsection{CM} \label{cm}
A Calabi-Yau threefold $X$ is said to have CM if the Mumford Tate group of the Hodge structure on $H^3(X,\QQ)$ is abelian (cf. \cite{Rohde-maximal}, section 6). As the Mumford-Tate group of $H^1(E,\QQ)$ is abelian, $X_f$ has $CM$ if and only if $H^1(C_f,\QQ)$ has an abelian Mumford-Tate group.
It is well-known that the Jacobian of $C_f$ is isogenous a product of two (isogenous)
elliptic curves. In fact,
the genus two curves $C_f$ also admit a Weierstrass equation of the form
$y^2=(x^3-a)(x^3-b)$ and scaling $x$ suitably one finds the equation
$y^2=(x^3-c)(x^3-c^{-1})$. This shows that $(x,y)\mapsto (1/x,y/x^3)$ is an involution of the curve and its eigenspaces in the tangent space of the Jacobian
define two elliptic curves in the Jacobian. Using the automorphism of order three 
one finds that these curves are isogenous. This implies that $H^1(C_f,\QQ)$ has an abelian Mumford Tate group if and only if these elliptic curves have CM in the classical sense, that is, that their endomorphism algebra is an order in an imaginary quadratic field. 
As these elliptic curves are also isogenous to the quotient of the genus two curve by the involution, 
this allows one to find explicit Calabi-Yau threefolds with CM.

\subsection{Hypergeometric differential equations} \label{hge}
It follows from Proposition \ref{cohX}
that the variation of Hodge structures given by the 
$H^3(X_f,\QQ)$ is the same as the variation of Hodge structures
of the $H^1(C_f,\QQ)$. This latter variation has been extensively studied (starting with Euler(!)) and we recall the main result.
(We found the results below in the unpublished PhD thesis by B.\ van der Marel, see also \cite{L}, section 4).
First of all, we assume that one of the four branch points of
the cyclic degree three cover $C_f\rightarrow \PP^1_t$ is the point $\infty\in\PP^1$. Then the curve $C_f$ is isomorphic to the curve
defined by
$$
C_\lambda\,:\; y^N\,=\,x^A(x-1)^B(x-\lambda)^C
$$
for some $\lambda\in \CC$ with $N=3$ and one can assume that $A=B=1$ and $C=2$. The holomorphic one forms on this curve are 
$\omega(0,0,0;1)$ and $\omega(0,0,1;2)$ where
$$
\omega(\alpha,\beta,\gamma;l)\;:=\;
\frac{x^\alpha(x-1)^\beta(x-\lambda)^\gamma}{y^l}\rmd x.
$$
To find the Picard-Fuchs equations, let:
$$
a\,:=\,-\alpha+(lA/N),\qquad
b\,:=\,-\beta+(lB/N),\qquad
c\,:=\,-\gamma+(lC/N).
$$
Then we have
$$
\omega(\alpha,\beta,\gamma;l)\;=\;
x^{-a}(x-1)^{-b}(x-\lambda)^{-c}\rmd x,\qquad
\frac{\partial}{\partial \lambda} \omega(\alpha,\beta,\gamma;l)\,=\,
\frac{c}{(x-\lambda)}\omega(\alpha,\beta,\gamma;l)
$$
etc.
An explicit computation shows that
$$
\left(\lambda(1-\lambda)\frac{\partial^2}{\partial \lambda^2}+
(a+c-\lambda(a+b+2c))\frac{\partial}{\partial \lambda}-c(a+b+c-1)\right)
\omega(\alpha,\beta,\gamma;l)\;=\;c\rmd h
$$
where $h$ is the rational function
$$
h\,:=\,\frac{x^{\alpha+1}(x-1)^{\beta+1}(x-\lambda)^{\gamma-1}}{y^l}\,=\,
x^{1-a}(x-1)^{1-b}(x-\lambda)^{-1-c}.
$$

As already explained in \cite{Rohde-maximal}, the fact that one obtains two degree two equations rather than one degree four equation follows from the fact that $\partial/\partial\lambda$ commutes with the automorphism $(x,y)\mapsto (x,\xi y)$ of $C_\lambda$. In particular, the monodromy of the Picard-Fuchs equation of the variation of Hodge structures
$\{H^1(C_\lambda,\ZZ)\}_\lambda$ is, up to conjugation, contained in the
subgroup of matrices with diagonal $2\times 2$ blocks in $Sp(4,\CC)$ and thus there is no point in $\PP^1$ with maximal unipotent monodromy.
It would be interesting to see if the family is the Mirror of a family of Calabi-Yau threefolds, these threefolds would have $h^{1,1}=1$ and $h^{2,1}=73$.

\

\end{document}